\documentclass[3p]{elsarticle}

\usepackage{amsmath}
\usepackage{amsthm}
\usepackage{amssymb}
\usepackage{bm}
\usepackage{accents}
\usepackage{tikz}
\usepackage{pgfplots}

\newtheorem{thm}{Theorem}

\newdefinition{rmk}{Remark}
\newproof{pf}{Proof}

\newcommand{\const}{\mathop{\rm const}\nolimits}

\numberwithin{equation}{section}

\graphicspath{{figs/}}

\journal{arXiv}

\begin{document}

\begin{frontmatter}

\title{Operator-difference schemes on non-uniform grids for second-order evolutionary equations}

\author{P.N. Vabishchevich\corref{cor1}\fnref{lab1,lab2}}
\ead{vab@ibrae.ac.ru}
\cortext[cor1]{Correspondibg author.}

\address[lab1]{Nuclear Safety Institute, Russian Academy of Sciences,
              52, B. Tulskaya, 115191 Moscow, Russia}
\address[lab2]{North-Eastern Federal University, 58, Belinskogo st, Yakutsk, 677000, Russia}

\begin{abstract}
The approximate solution of the Cauchy problem for second-order evolution equations is performed, first of all, using three-level time approximations.
Such approximations are easily constructed and relatively uncomplicated to investigate when using uniform time grids.
When solving applied problems numerically, we should focus on approximations with variable time steps.
When using multilevel schemes on non-uniform grids, we should maintain accuracy by choosing appropriate approximations and ensuring the approximate solution's stability.
In this paper, we construct unconditionally stable first- and second-order accuracy schemes on a non-uniform time grid for the approximate solution of the Cauchy problem for a second-order evolutionary equation.
We use a special transformation of the original second-order differential-operator equation to a system of first-order equations.
For the system of first-order equations, we apply standard two-level time approximations.
We obtained stability estimates for the initial data and the right-hand side in finite-dimensional Hilbert space.
Eliminating auxiliary variables leads to three-level schemes for the initial second-order evolution equation.
Numerical experiments were performed for the test problem for a one-dimensional in space bi-parabolic equation.
The accuracy and stability properties of the constructed schemes are demonstrated on non-uniform grids with randomly varying grid steps.
\end{abstract}


\end{frontmatter}

\section{Introduction}\label{sec:1}

Two-level schemes (Runge-Kutta methods) \cite{butcher2008numerical,hairer1993solving} are widely used for the approximate solution of the Cauchy problem for systems of ordinary differential equations.
The second class of methods is based on multilevel time approximations (multistep methods) when transitioning to a new time level involves a known solution on two or more preceding time levels.
Two- and three-level approximations are most common when numerically solving nonstationary boundary value problems for partial differential equations.
When considering the Cauchy problem for second-order evolutionary equations, it is natural to focus on three-level approximations.

Multilevel schemes on uniform time grids are the easiest to construct and investigate.
The computational practice focuses on using non-uniform grids for accounting for local features of the solution by choosing the time step.
When using multilevel approximations on non-uniform grids, there are two main problems.
The first one is related to the choice of approximations in time to ensure the required accuracy of the approximated solution.
The second problem is generated by ensuring the conditions for the stability of the solution under the selected approximations.

The results of the theory of stability (correctness) of operator-difference schemes \cite{Samarskii1989} are widely used for the construction and investigation of time-difference approximations in the numerical solution of applied problems.
In this case, the difference scheme is considered in finite-dimensional Hilbert spaces; the scheme is written in a unified (canonical) form for which stability conditions exist.
The necessary and sufficient conditions for stability of two- and three-level schemes on uniform time grids have been obtained (see, for example, \cite{SamarskiiMatusVabischevich2002}).

For two-level schemes, the transition to non-uniform grids is editorial.
In the case of three-level schemes, the situation is much more complicated, both in the choice of approximations and in the study of stability.
Let us note some available results in this direction in the particular case when the time step change is rare.
When considering the stability of three-level difference schemes, we must consider the problem separately when the time step increases and when it decreases.
We have found (see, e.g., \cite{samarskii2001stability,matus2001three,cho2013stability}) that when the time step changes, the a priori stability estimates become worse, with the constant in the a priori estimate increasing when the time step changes in proportion to the time step ratio.

We can proceed to a system of first-order evolution equations for high-order evolution equations.
We can use the usual two-level approximations on non-uniform grids for the resulting vector problem.
The specificity of the problem is considered by choosing auxiliary unknowns so that the two-level schemes for the system of first-order equations have the necessary properties of accuracy and stability.
This paper implements such an approach for the Cauchy problem for a general second-order differential-operator equation with self-adjoint constant positively determined operators in a finite-dimensional Hilbert space.
Unconditionally stable first- and second-order approximation schemes for the introduced system of first-order evolution equations have been constructed.
Algebraically equivalent unconditionally stable three-level schemes of the second order are obtained when the auxiliary variable is eliminated.

We give a brief content of the paper.
Section 2 sets up a Cauchy problem for a second-order evolutionary equation.
A typical result on the stability of a three-level difference scheme on a uniform time grid is established.
The original second-order equation is written down in Section 3 as a system of first-order equations when a new unknown quantity is introduced.
An a priori estimate for the solution of the resulting vector Cauchy problem is proved.
The critical part of the paper is Section 3.
The unconditionally stable two-level schemes for the system of equations are constructed.
With the elimination of auxiliary variables, a three-level scheme on a non-uniform time grid was obtained for the original second-order evolution equation.
Section 5 presents the results of numerical experiments on a uniform and non-uniform grids.
The Cauchy problem for the one-dimensional in space bi-parabolic equation is considered a test problem.
We summarize the results of our study in Section 6.

\section{Problem formulation}\label{sec:2}

We consider the Cauchy problem for a second-order evolutionary equation in the real finite-dimensional Hilbert space $H$.
The function $u(t)$ satisfies the equation
\begin{equation}\label{2.1}
C \frac{d^2 u}{d t^2} + B\frac{d u}{d t} + A u = f(t),
\quad t > 0 ,
\end{equation}
and the initial conditions
\begin{equation}\label{2.2}
u(0) = u_0 ,
\quad \frac{d u}{d t}(0) = u'_0 .
\end{equation}
The linear constants (independent of $t$) operators $A, B, C$ are self-adjoint and positively definite:
\begin{equation}\label{2.3}
A = A^* \geq \nu_A I, \quad \nu_A > 0 ,
\quad B = B^* \geq \nu_B I, \quad \nu_B > 0,
\quad C = C^* \geq \nu_C I, \quad \nu_C > 0,
\end{equation}
where $I$ is a unit operator in $H$.
We will use the usual notations $(\cdot, \cdot)$ and $\|\cdot \|$ for the scalar product and norm in $H$.
For a self-adjoint positively defined operator $D$, we define a Hilbert space $H_D$ with scalar product and norm
$(u, v)_D = (D u, v), \ \| u \|_D = (u, v)_D^{1/2}.$

We get such problems, in particular, when using the finite element method or the finite volume method for approximation by space in the numerical study of nonstationary boundary value problems for partial differential equations.
We can associate equation (\ref{2.1}) at $C = 0$ with a parabolic equation and a hyperbolic equation at $B = 0$.

Let us give a simple a priori estimate for the solution of the problem (\ref{2.1})---(\ref{2.3}); this estimate will guide us in the study of time approximations.
By multiplying equation (\ref{2.1}) scalarly in $H$ by $d u / dt$, we obtain the equality
\[
\frac{1}{2} \frac{d}{d t} \Big ( \Big \| \frac{d u}{d t} \Big \|_C^2 + \|u\\|A^2 \Big ) + \Big \| \frac{d u}{d t} \Big \|_B^2 =
\Big ( f, \frac{d u}{d t} \Big ) .
\]
For the right-hand side, we use the estimate
\[
\Big ( f, \frac{d u}{d t} \Big ) \leq \Big \| \frac{d u}{d t} \Big \|_B^2 + \frac{1}{4} \|f\\|2_{B^{-1}} .
\]
Given this inequality, Gronwall's lemma gives us an a priori estimate
\begin{equation}\label{2.4}
\Big \| \frac{d u}{d t} (t) \Big \|_C^2 + \|u(t)\|_A^2 \leq \| v_0 \|_C^2 + \|u_0\|_A^2
+ \frac{1}{2} \int_{0}^{t} \|f(\theta)\|^2_{B^{-1}} d \theta .
\end{equation}

In the approximate solution of the Cauchy problem (\ref{2.1}), (\ref{2.2}), it is natural to focus on unconditionally stable three-level difference schemes.
It is easiest to work with uniform time grids when the time step $\tau = \const$ and let $y^n=y(t_n), \ t_n = n\tau$, $n =0, 1, \ldots$.
A three-level difference scheme with weight
$\sigma = \mathrm{const} > 0$ has the form
\begin{equation}\label{2.5}
C \frac{ y^{n+1} - 2 y^{n} + y^{n-1}}{\tau^2 } + B \frac{ y^{n+1} - y^{n-1}}{2\tau } + A y^{n+\sigma} =
f^{n}
\quad n = 1,2,\ldots,
\end{equation}
\begin{equation}\label{2.6}
y^0 = u^0 ,
\quad y^1 = \accentset{\sim }{ u}^1 ,
\end{equation}
when using the notations
\[
y^{n+\sigma} = \sigma y^{n+1} + (1-2\sigma ) y^{n} + \sigma y^{n-1} .
\]
Given the initial conditions (\ref{2.2}), we put
\[
u^{0} = u_0,
\quad \accentset{\sim }{ u}^1 = u_0 + \tau u'_0
\quad (\accentset{\sim }{ u}^1 \approx u^1(\tau) \approx u_0 + \tau u'_0 ).
\]

The three-level scheme (\ref{2.5}), (\ref{2.6}) has the second order of approximation by $\tau$.
The theory of stability of operator-difference schemes gives the following conditions for unconditional stability of this scheme: $\sigma \geq 0.25$.

\begin{thm}\label{t-1}
The three-level difference scheme (\ref{2.3}), (\ref{2.5}), (\ref{2.6})
is unconditionally stable for $\sigma \geq 0.25$.
Under these constraints, for the approximate solution of (\ref{2.1}), (\ref{2.2}), the a priori estimate is
\begin{equation}\label{2.7}
\begin{split}
\Big \|\frac{y^{n+1} - y^{n}}{\tau } \Big \|_D^2 & +
\Big \|\frac{ y^{n+1} + y^{n}}{2} \Big \|^2_{ A} \\
& \leq \Big \|\frac{\accentset{\sim }{u}^1 - u^0}{\tau } \Big \|_D^2 +
\Big \|\frac{\accentset{\sim }{u}^1 + u^0}{2} \Big \|^2_{ A}
+ \frac{1}{2} \sum_{k=1}^{n}\tau \|f^{k}\|^2_{B^{-1}} ,
\quad n = 1,2,\ldots ,
\end{split}
\end{equation}
where
\[
D = C + \Big (\sigma - \frac{1}{4} \Big ) \tau^2 A .
\]
\end{thm}

\begin{pf}
Given the equality
\[
\sigma y^{n+1} + (1-2\sigma ) y^{n} + \sigma y^{n-1} =
\frac{1}{4} (y^{n+1} + 2 y^{n} + y^{n-1}) + \Big (\sigma - \frac{1}{4} \Big ) (y^{n+1} - 2 y^{n} + y^{n-1} ) ,
\]
we write equation (\ref{2.5}) as
\begin{equation}\label{2.8}
D \frac{ y^{n+1} - 2 y^{n} + y^{n-1}}{\tau^2 } + B \frac{ y^{n+1} - y^{n-1}}{2\tau } + \frac{1}{4} A (y^{n+1} + 2 y^{n} + y^{n-1}) =
f^{n} .
\end{equation}
Let's denote by
\[
s^n = y^n - y^{n-1},
\quad r^n = \frac12 (y^n + y^{n-1}),
\]
then
\[
\begin{split}
y^{n+1} - y^{n-1} & = s^{n+1} + s^n , \\
y^{n+1} - 2 y^{n} + y^{n-1} & = s^{n+1} - s^n , \\
y^{n+1} + 2 y^{n} + y^{n-1} & = 2 (r^{n+1} + r^n) .
\end{split}
\]
Equation (\ref{2.8}) takes the form
\begin{equation}\label{2.9}
D \frac{s^{n+1} - s^n}{\tau^2 } + B \frac{s^{n+1} + s^n}{2\tau } + A \frac{r^{n+1} + r^n}{2} =
f^{n} .
\end{equation}
Let's multiply equation (\ref{2.9}) by $s^{n+1} + s^n = 2(r^{n+1} - r^{n})$; this give
\[
\frac{1}{\tau^2 } \big (\|s^{n+1}\|_D^2 - \|s^{n}\|_D^2 \big ) + \|r^{n+1} \|^2_{ A} - \|r^{n} \|^2_{ A}
+ \frac{1}{2 \tau } \|s^{n+1} + s^n\|_B^2 = (f^n, s^{n+1} + s^n) .
\]
Using the introduced notations and the inequality
\[
(f^n, s^{n+1} + s^n) \leq \frac{1}{2 \tau } \|s^{n+1} + s^n\|_B^2 + \frac{\tau}{2 } \|f^n\|^2_{B^{-1}} ,
\]
we get an estimate of the solution on one level in time
\begin{equation}\label{2.10}
\Big \|\frac{y^{n+1} - y^{n}}{\tau } \Big \|_D^2 +
\Big \|\frac{y^{n+1} + y^{n}}{2} \Big \|^2_{ A}
\leq \Big \|\frac{y^{n} - y^{n-1}}{\tau } \Big \|_D^2 + \Big \|\frac{y^{n} + y^{n-1}}{2} \Big \|^2_{ A}
+ \frac{1}{2} \tau \|f^n\|^2_{B^{-1}} .
\end{equation}
We obtain the provable estimate (\ref{2.7}) by applying the difference analog of Gronwall's lemma.
\end{pf}

\begin{rmk}
We can consider the a priori estimate (\ref{2.7}) as a direct discrete analogue of the estimate (\ref{2.4}) for solving the Cauchy problem (\ref{2.1})--(\ref{2.2}).
\end{rmk}

In computational practice, in simulating nonstationary processes, we should focus on using non-uniform time grids.
For multilevel schemes on non-uniform grids, there are problems with maintaining the accuracy of the approximations.
The stability conditions of schemes on non-uniform grids also deteriorate.
In particular, when considering three-level operator-difference schemes, we have stability estimates that depend on the time step ratio and how often the time step change occurs.
Our goal is to construct nonuniform time approximations for nonstationary problems (\ref{2.1})--(\ref{2.3}), which would give unconditionally stable schemes of the second order of accuracy.

\section{System of first-order equations}\label{sec:3}

Constructing approximations on non-uniform grids for first-order evolution equations is solved using various variants of two-level schemes (Runge-Kutta methods).
On the other hand, we can always proceed from evolution equations of the second and higher order to a system of evolution equations of the first order.
In this case, the focus should be on choosing new variables.
We consider such possibilities of constructing approximations on non-uniform time grids for the approximate solution of the problem (\ref{2.1})--(\ref{2.3}).

The simplest version of introducing a new unknown quantity $v = du/dt$ for equation (\ref{2.1}), (\ref{2.3}) is not of great interest.
We put
\[
v = C \frac{d u}{d t} + B u
\]
and write equation (\ref{2.1}) as a system of equations
\begin{equation}\label{3.1}
C \frac{d u}{d t} + B u - v = 0 ,
\end{equation}
\begin{equation}\label{3.2}
A^{-1} \frac{d v}{d t} + u = A^{-1} f ,
\quad t > 0 .
\end{equation}
For the system (\ref{3.1}), (\ref{3.2}), a Cauchy problem is posed when the initial conditions are
\begin{equation}\label{3.3}
u(0) = u_0,
\quad v(0) = v_0 ,
\end{equation}
where given (\ref{2.2}) we have $v_0 = C u'_0 + B u_0$.

It is convenient for us to write the system of equations (\ref{3.1}), (\ref{3.2}) as one equation for vector quantities.
Define vectors $\bm u = \{u, v\}$ and $\bm f = \{0, A^{-1} f\}$ and from (\ref{3.1})--(\ref{3.2}) we obtain the Cauchy problem
\begin{equation}\label{3.4}
\bm B \frac{d \bm u}{d t} + \bm A \bm u = \bm f ,
\quad t > 0 ,
\end{equation}
\begin{equation}\label{3.5}
\bm u(0) = \bm u_0 ,
\end{equation}
where $\bm u_0 = \{u_0, v_0 \}$.
For the operator matrices $\bm B$ and $\bm A$, we obtain the expressions
\[
\bm B = \mathrm{diag} \, \big (C, A^{-1} \big),
\quad \bm A = \left (\begin{array}{cc}
B & - I \\
I & 0 \\
\end{array}
\right ) .
\]

The problem (\ref{3.4}), (\ref{3.5}) we will consider on the direct sum of spaces $\bm H = H \oplus H$. For $\bm u_1, \bm u_2 \in \bm H$ the scalar product and norm are defined by the expressions
\[
 (\bm u_1, \bm u_2) = (u_1, v_1) + (u_2, v_2),
 \quad \|\bm u\| = (\bm u, \bm u)^{1/2} .
\]
Given (\ref{2.3}), we obtain
\begin{equation}\label{3.6}
 \bm B = \bm B^* > 0,
 \quad \bm A \geq 0 . 
\end{equation} 
Such good properties of the $\bm B$ and $\bm A$ operators are provided by choosing the second equation of the system as (\ref{3.2}).

To obtain an a priori estimate for the solution of (\ref{3.4}), (\ref{3.5}), we multiply the equation (\ref{3.4}) by $\bm u$ scalarly in $\bm H$.
Considering (\ref{3.6}), this gives 
\[
 \|\bm u\|_{\bm B} \frac{d}{d t}  \|\bm u\|_{\bm B} \leq (\bm f , \bm u) .  
\]
Taking into account of
\[
 (\bm f , \bm u) \leq \|\bm \varphi \|_{\bm B^{-1}} \|\bm u\|_{\bm B} ,
\]  
we get
\begin{equation}\label{3.7}
 \|\bm u(t)\|_{\bm B} \leq \|\bm u_0\|_{\bm B} + \int_{0}^{t} \|\bm f(s)\|_{\bm B^{-1}} d s .
\end{equation} 
In our case
\[
 \|\bm u(t)\|_{\bm B} = \Big (\|u(t)\|_C^2 + \|v(t)\|_{A^{-1}}^2 \Big)^{1/2} , 
 \quad \|\bm v_0\|_{\bm B} = \Big (\|u_0\|_C^2 + \|v_0\|_{A^{-1}}^2 \Big)^{1/2} ,
 \quad \|\bm f(t)\|_{\bm B^{-1}} = \|f(t)\|_{A^{-1}} ,
\] 
so for the problem (\ref{3.1})--(\ref{3.3}) we have an estimate
\begin{equation}\label{3.8}
 \Big (\|u(t)\|_C^2 + \|v(t)\|_{A^{-1}}^2 \Big)^{1/2} \leq \Big (\|u_0\|_C^2 + \|v_0\|_{A^{-1}}^2 \Big)^{1/2} +
 \int_{0}^{t} \|f(s)\|_{A^{-1}} d s .
\end{equation} 

\begin{rmk}
For the primary problem (\ref{2.1})--(\ref{2.3}), the estimate (\ref{3.8}) is established by scalarly multiplying equation (\ref{2.1}) by
\[
A^{-1} v = A^{-1} \Big ( C \frac{d u}{d t} + B u \Big ) .
\]
\end{rmk}

\begin{rmk}
The positive definiteness condition of the operator $B$ (see (\ref{2.3})) can be relaxed.
The estimates (\ref{3.7}), (\ref{3.8}) take place even when $\nu_B = 0$.
In this case, we must replace the estimate of (\ref{2.4}) for the solution of (\ref{2.1}), (\ref{2.2}) by (\ref{3.8}).
\end{rmk}

\section{Operator-difference schemes on non-uniform grids}\label{sec:4}

Two-level schemes are most easily constructed to approximate the problem (\ref{3.4}), (\ref{3.5}).
We consider approximations on a non-uniform time grid when
\[
t_{n+1} = t_{n} + \tau_{n+1},
\quad n = 0,1, \ldots,
\quad t_0 = 0 .
\]
For the approximate solution, we will use a two-level scheme with weight $\sigma = \const$
\begin{equation}\label{4.1}
\bm B \frac{\bm y^{n+1} - \bm y^{n}}{\tau_{n+1}} + \bm A \bm y^{n+\sigma} = \bm f^{n+\sigma}
\quad n = 0,1,\ldots,
\end{equation}
\begin{equation}\label{4.2}
\bm y^0 = \bm u_0 ,
\end{equation}
when using the notations
\[
\bm y^{n} = (y^n, w^n),
\quad \bm y^{n+\sigma} = \sigma \bm y^{n+1} + (1-\sigma ) \bm y^{n} .
\]
For a sufficiently smooth solution of $\bm u(t)$ of the problem (\ref{3.4}), (\ref{3.5}), the difference scheme (\ref{4.1}), (\ref{4.2}) has first order approximation by $\tau$ for $\sigma \neq 0.5$.
For $\sigma = 0.5$ we have a Crank-Nicolson scheme; this scheme approximates (\ref{3.4}), (\ref{3.5}) with second order on $\tau$.
Necessary and sufficient stability conditions for such schemes are well known in the theory of stability (correctness) of operator-difference schemes.
The following statement guarantees sufficient conditions for unconditional stability.

\begin{thm}\label{t-2}
At $\sigma \geq 0.5$ the two-level operator-difference scheme (\ref{3.4}), (\ref{3.5}) is unconditionally stable.
In this case, there is an a priori estimate for the approximate solution
\begin{equation}\label{4.3}
\|\bm y^{n+1}\|_{\bm B} \leq \|\bm u_0\|_{\bm B} + \sum_{k=0}^{n}\tau_{n+1} \|\bm f^{k+\sigma}\|_{\bm B^{-1}} ,
\quad n = 0,1,\ldots .
\end{equation}
\end{thm}

\begin{pf}
We will multiply equation (\ref{4.1}) scalarly by $\tau \bm y^{n+\sigma}$. Given the non-negativity of the operator $\bm A$, we have the inequality
\begin{equation}\label{4.4}
(\bm B (\bm y^{n+1} - \bm y^{n}), \bm y^{n+\sigma}) \leq \tau \|\bm f^{n+\sigma}\|_{\bm B^{-1}} \|\bm y^{n+\sigma }\|_{\bm B} .
\end{equation}
At $\sigma \geq 0.5$ we use (see lemma 1 in \cite{vabishchevich2013flux}) for the left-hand side of the inequality (\ref{4.4}) the estimate
\[
\big (\bm B (\bm y^{n+1} - \bm y^{n}), y^{n+\sigma} \big ) \geq
\big (\|\bm y^{n+1}\|_{\bm B} - \|\bm y^{n}\|_{\bm B} \big ) \|\bm y^{n+\sigma}\|_{\bm B} .
\]
By doing so, we get an estimate for the solution on the new time level
\[
\|\bm y^{n+1}\|_{\bm B} \leq \|\bm y^{n}\|_{\bm B} + \tau \|\bm f^{n+\sigma}\|_{\bm B^{-1}} .
\]
The difference analog of Gronwall's lemma gives us an a priori estimate (\ref{4.3}), which is the grid analog of the estimate (\ref{3.7}).
\end{pf}

Applying the vector scheme (\ref{3.4}), (\ref{3.5}) for the approximate solution of the problem (\ref{3.1})--(\ref{3.3}) leads us to the operator-difference scheme
\begin{equation}\label{4.5}
C \frac{y^{n+1} - y^{n}}{\tau_{n+1}} + B (\sigma y^{n+1} + (1-\sigma) y^{n}) - (\sigma w^{n+1} + (1-\sigma) w^{n}) = 0,
\end{equation}
\begin{equation}\label{4.6}
\frac{w^{n+1} - w^{n}}{\tau_{n+1}} + A (\sigma y^{n+1} + (1-\sigma) y^{n}) = f^{n+\sigma} ,
\end{equation}
\begin{equation}\label{4.7}
y^0 = u_0,
\quad w^0 = v_0 .
\end{equation}
Under these conditions, the stability estimate (\ref{4.3}) takes the form
\begin{equation}\label{4.8}
\Big (\|y^{n+1} \|_C^2 + \|w^{n+1} \|_{A^{-1}}^2 \Big)^{1/2} \leq \Big (\|u_0\||_C^2 + \|v_0\|_{A^{-1}}^2 \Big)^{1/2} +
\sum_{k=0}^{n}\tau_{k+1} \|f^{k+\sigma}\|_{A^{-1}} .
\end{equation}
The estimate (\ref{4.8}) is the discrete analog of the a priori estimate (\ref{3.8}) for the problem (\ref{3.1})--(\ref{3.3}).

The computational realization of the transition to a new level in time for the system of equations (\ref{3.1}), (\ref{3.2}) is provided as follows.
From equation (\ref{4.6}), we get the representation
\begin{equation}\label{4.9}
w^{n+1} = - \sigma \tau_{n+1} A y^{n+1} + \chi_1^n,
\end{equation}
at a known
\[
\chi_1^n = w^n - (1-\sigma)\tau_{n+1} A y^{n} + \tau_{n+1} f^{n+\sigma} .
\]
Similarly, from equation (\ref{4.5}) we have
\begin{equation}\label{4.10}
(C + \sigma \tau_{n+1}B )y^{n+1} = \sigma \tau_{n+1} w^{n+1} + \chi_2^n,
\end{equation}
where
\[
\chi_2^n = (C - (1-\sigma)\tau_{n+1} B) y^{n} + (1-\sigma)\tau_{n+1} w^n .
\]
Substituting (\ref{4.9}) into (\ref{4.10}) gives the equation to determine $y^{n+1}$
\begin{equation}\label{4.11}
(C + \sigma \tau_{n+1}B + \sigma^2 \tau^2_{n+1} A) y^{n+1} = \chi^n,
\end{equation}
at the right-hand side
\[
\chi^n = \sigma \tau_{n+1} \chi_1^n + \chi_2^n.
\]
After solving the equation (\ref{4.11}), the expression (\ref{4.9}) is used to determine $w^{n+1}$.

When applying the standard three-level scheme on a uniform grid (\ref{2.5}), (\ref{2.6}) while determining $y^{n+1}$ we solve the problem
\[
\Big (C + \frac{1}{2} \tau B + \sigma \tau^2 A \Big ) y^{n+1} = \widetilde{\chi}^n,
\]
with the right-hand side
\[
\widetilde{\chi}^n = C (2 y^{n} - y^{n-1}) + \frac{1}{2} \tau B y^{n-1} - \tau^2 A ((1-2\sigma) y^n + \sigma y^{n-1}) + \tau^2 f^n.
\]
A comparison with the problem (\ref{4.11}) at $\tau_{n+1} = \tau$ shows that these problems are close (left-hand sides are the same) if we put $\sigma = 0.25$ in equation (\ref{2.5}) and set $\sigma = 0.5$ in the system (\ref{4.5}), (\ref{4.6}).

The two-level scheme (\ref{4.5}), (\ref{4.6}) for two variables $y$ and $w$ can be written as a three-level scheme only for $y$.
To do this, we write the equation (\ref{4.5}) on two adjacent time levels and write their difference in the form
\begin{equation}\label{4.12}
\begin{split}
C \Big (\frac{y^{n+1} - y^{n}}{\tau_{n+1}} - \frac{y^{n} - y^{n-1}}{\tau_{n}} \Big )
& + B \big (\sigma (y^{n+1}-y^{n}) + (1-\sigma) (y^{n}- y^{n-1}) \big ) \\
& = \sigma (w^{n+1}-w^{n}) + (1-\sigma) (w^{n}- w^{n-1}) .
\end{split}
\end{equation}
For the right-hand side, the equation (\ref{4.6}) is also used at two levels in time; this gives
\begin{equation}\label{4.13}
\begin{split}
\sigma (w^{n+1}-w^{n}) & + (1-\sigma) (w^{n}- w^{n-1}) \\
& = - \sigma \tau_{n+1} A (\sigma y^{n+1} + (1-\sigma) y^{n}) - (1-\sigma) \tau_{n} A (\sigma y^{n} + (1-\sigma) y^{n-1}) \\
& + \sigma \tau_{n+1} f^{n+\sigma} + (1-\sigma) \tau_{n} f^{n-1+\sigma} .
\end{split}
\end{equation}
From (\ref{4.12}), (\ref{4.13}), we get a three-level scheme for the variable $y$.

In the most interesting case $\sigma = 0.5$, when the scheme (\ref{4.5}), (\ref{4.6}) has the second order of accuracy, we obtain
\begin{equation}\label{4.14}
\begin{split}
C \frac{2}{\tau_{n}+\tau_{n+1}} \Big (\frac{y^{n+1} - y^{n}}{\tau_{n+1}} - \frac{y^{n} - y^{n-1}}{\tau_{n}} \Big )
& + B \Big ( \overline{\sigma}_{n+1} \frac{y^{n+1} - y^{n}}{\tau_{n+1}} + (1- \overline{\sigma}_{n+1}) \frac{y^{n} - y^{n-1}}{\tau_{n}} \Big ) \\
& + A \Big ( \overline{\sigma}_{n+1} \frac{y^{n+1} + y^{n}}{2} + (1- \overline{\sigma}_{n+1}) \frac{y^{n} + y^{n-1}}{2} \Big ) \\
& = \overline{\sigma}_{n+1} f^{n+1/2} + (1-\overline{\sigma})_{n+1} f^{n-1/2} .
\end{split}
\end{equation}
The variable weight $\overline{\sigma}_{n+1}$ in (\ref{4.14}) is determined by the expression
\[
\overline{\sigma}_{n+1} = \frac{\tau_{n+1}}{\tau_{n}+\tau_{n+1}} .
\]
Given that for a sufficiently smooth function $\varphi (t)$
\[
\overline{\sigma}_{n+1} \varphi(t^{n+1/2}) + (1-\overline{\sigma})_{n+1} \varphi(t^{n-1/2}) =
\varphi(\overline{t}^{n}) + \mathcal{O} \big((\tau_{n+1}+\tau_{n})^2 \big ) ,
\quad \overline{t}^{n} = \frac{1}{2} (t^{n+1} + t^{n-1}) ,
\]
we relate the scheme (\ref{4.14}) to an approximation of equation (\ref{2.1}) on an non-uniform grid at $t = \overline{t}^{n}$.

\section{Numerical experiments}\label{sec:5}

Among the generalized heat conduction models \cite{joseph1989heat}, we will distinguish the bi-parabolic equation \cite{fushchich1990new}.
In this case, the homogeneous second-order evolutionary equation has the form
\begin{equation}\label{5.1}
\Big (\frac{d}{d t} + D \Big ) u + \alpha \Big (\frac{d}{d t} + D \Big )^2 u = 0,
\quad t > 0 ,
\end{equation}
where
\[
D = D^* \geq \nu_D I,
\quad \nu_D > 0,
\quad \alpha = \const > 0 .
\]
We take the initial conditions for equation (\ref{5.1}) as
\begin{equation}\label{5.2}
u(0) = u_0,
\quad \frac{d u}{d t} (0) = 0 .
\end{equation}
Equation (\ref{5.1}) is written as a second-order evolutionary equation (\ref{2.1}) at
\[
C = \alpha I,
\quad B = I + 2 \alpha D,
\quad A = D + \alpha D^2 .
\]

The solution of the problem (\ref{5.1}), (\ref{5.2}) is
\begin{equation}\label{5.3}
u(t) = \Big ( I + \alpha \Big (1 - \exp \Big (- \frac{t}{\alpha}\Big ) D \Big ) \exp (- D t) u_0 .
\end{equation}
It follows from the representation (\ref{5.3}) that the relaxation of the initial state by the exponential law is provided not only by the positively defined operator $D$ but also by the square of the parabolic operator (parameter $\alpha$).

The principal goal of our study is to construct time approximations.
We can trace the key points when considering one-dimensional space problems.
It is optional to move to more interesting multidimensional problems.
We use sufficiently detailed computational grids over space when we can neglect computational errors of approximations over space.

Let us assume that in dimensionless variables $x \in [0,1]$ and introduce a uniform grid with step $h$.
Let us denote by $\omega$ the set of internal grid nodes in space.
The test problem is considered under homogeneous boundary conditions of the first kind.
The difference operator $- D$ is associated with the second derivative operator on $x$.
For $y(x) = 0, \ x \notin \omega$, we define the operator $D$ by the relation
\[
D y = - \frac{1}{h^2} \big (y(x+h) - 2 y(x) + y(x-h) \big ),
\quad x \in \omega .
\]
To be able to compare the calculation results with data from other generalized thermal conductivity models (see \cite{vabishchevich2022numerical}), we restrict ourselves to the initial conditions
\[
u_0(x) = \left \{ \begin{array}{rr}
x, & 0 < x \leq 0.5 , \\\
0, & 0.5 < x < 1 , \\
\end{array}
\right .
\quad x \in \omega .
\]
We performed calculations on a grid over space with $h=2\cdot 10^{-3}$ for $0 < t \leq T, \ T = 0.1$.
We used the scheme (\ref{4.5}), (\ref{4.6}) with $\sigma = 0.5$.
For the second initial condition (\ref{4.7}) taking into account of (\ref{5.2}), we have $v_0 = u_0 + 2\alpha D u_0$.

The new level problem (\ref{4.11}) for finding $y^{n+1}$ under our conditions is
\[
R y^{n+1} = \chi^n,
\]
where
\[
R = \alpha I + \frac{\tau_{n+1}}{2} (I + 2 \alpha D) + \frac{\tau^2_{n+1}}{4} (D + \alpha D^2) .
\]
For the operator $R$, we use the factorization
\[
R = (p D + I) \big ( \alpha p D + (p + \alpha) I \big ),
\quad p = \frac{\tau_{n+1}}{2} .
\]
To determine $y^{n+1}$, two common discrete problems with operators $D + cI$ at corresponding values of $c = \const > 0$ are solved.

\begin{figure}
\centering
\includegraphics[width=0.75\linewidth]{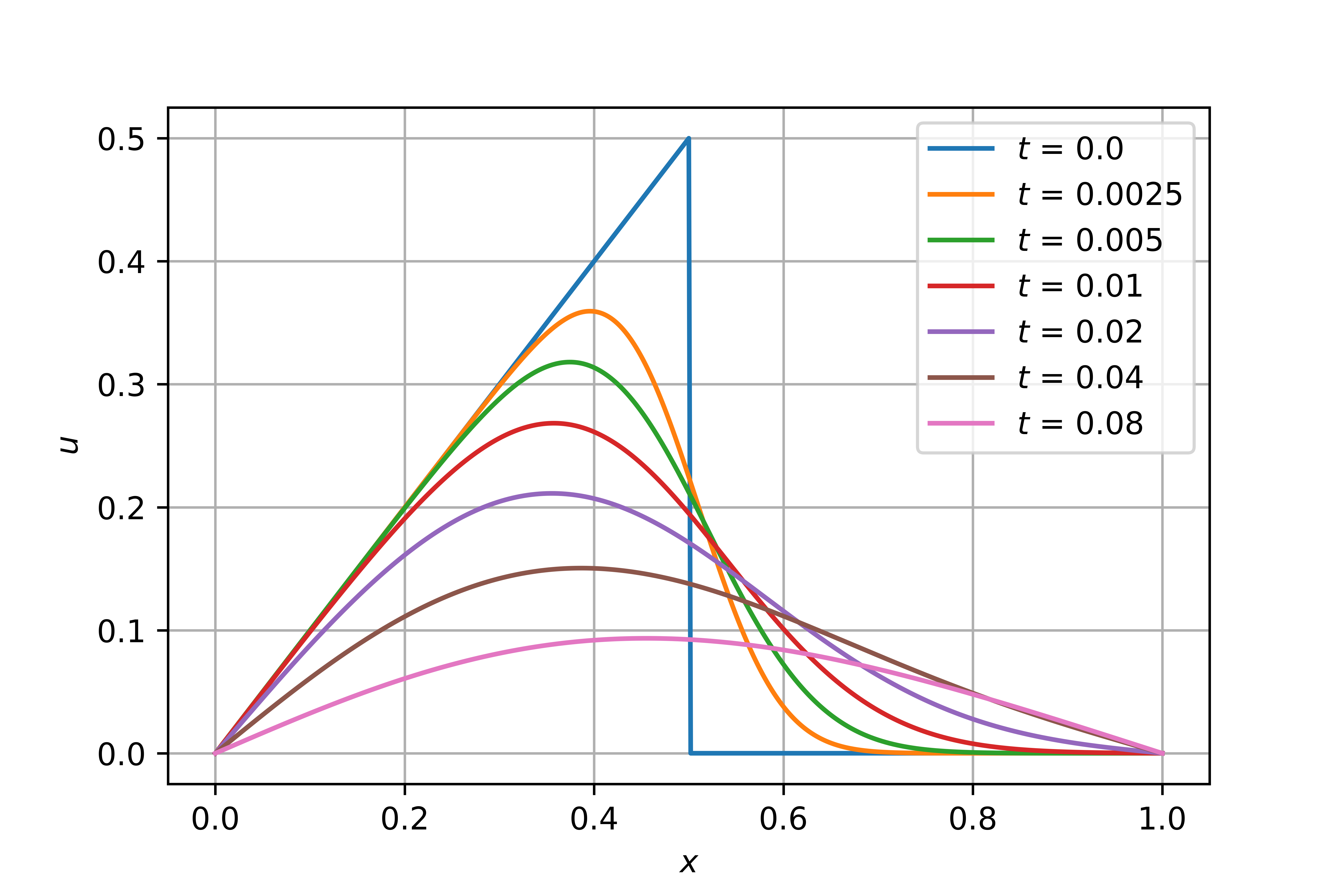} 
\caption{Exact solution of the problem at different points  in time when  $\alpha = 0$.}
\label{f-1}
\end{figure}

\begin{figure}
\centering
\includegraphics[width=0.75\linewidth]{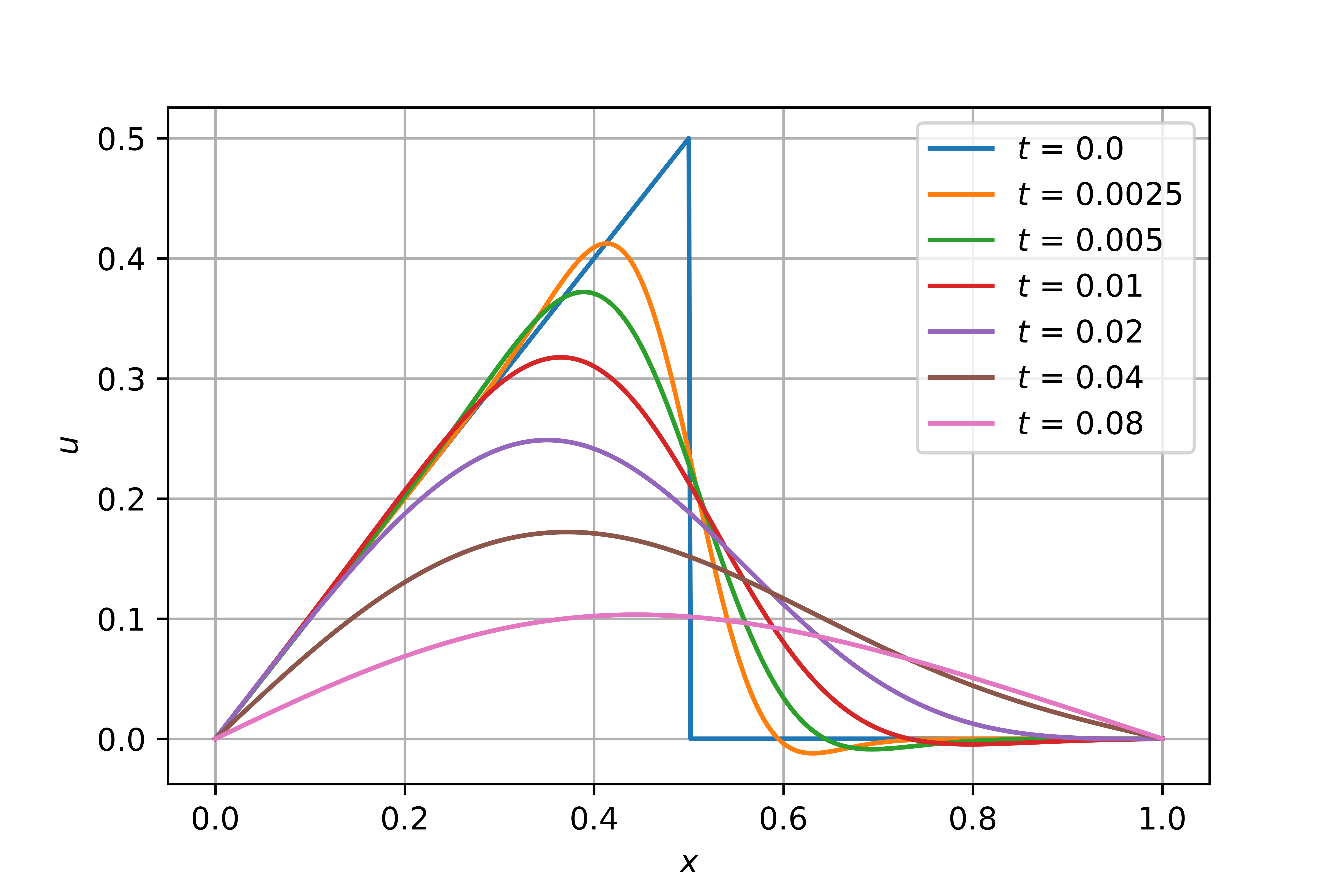} 
\caption{Exact solution of the problem at different points  in time when  $\alpha = 0.01$.}
\label{f-2}
\end{figure}

\begin{figure}
\centering
\includegraphics[width=0.75\linewidth]{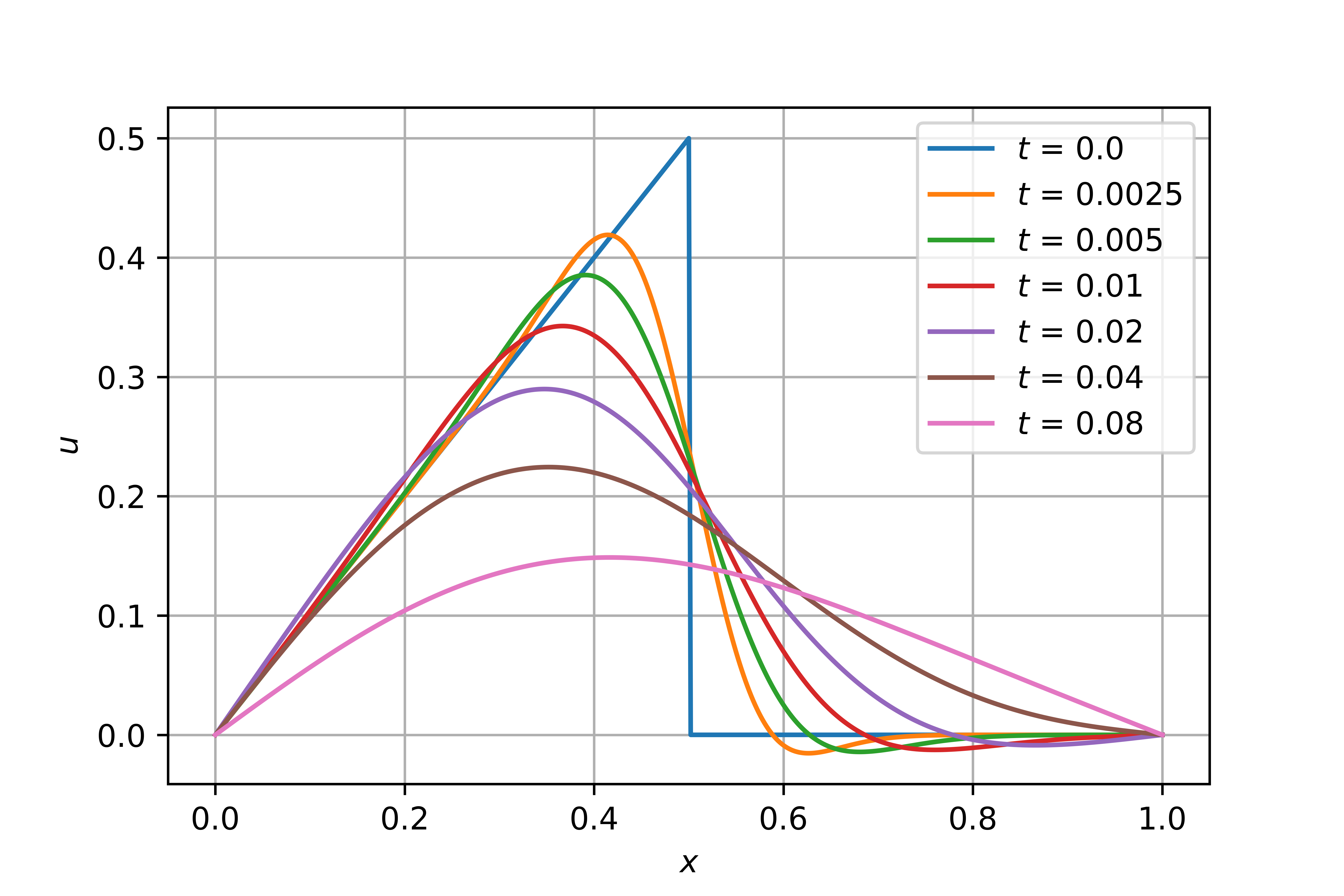} 
\caption{Exact solution of the problem at different points  in time when  $\alpha = 0.1$.}
\label{f-3}
\end{figure}

The exact solution of the problem (\ref{5.1}), (\ref{5.2}) on the used by us grid over space at different values of the parameter $\alpha$ is shown in Figures \ref{f-1}--\ref{f-3}.
The solution of the standard parabolic problem ($\alpha=0$, Fig.\ref{f-1}) is characterized by a monotonic decrease over time of the positive maximum of the solution, which is a consequence of the maximum principle.
The most striking manifestation of the hyperbolicity of the problem at positive $\alpha$ is (see Fig.\ref{f-2} and Fig.\ref{f-3}) the appearance of a region of negative solutions.

\begin{figure}
\centering
\includegraphics[width=0.75\linewidth]{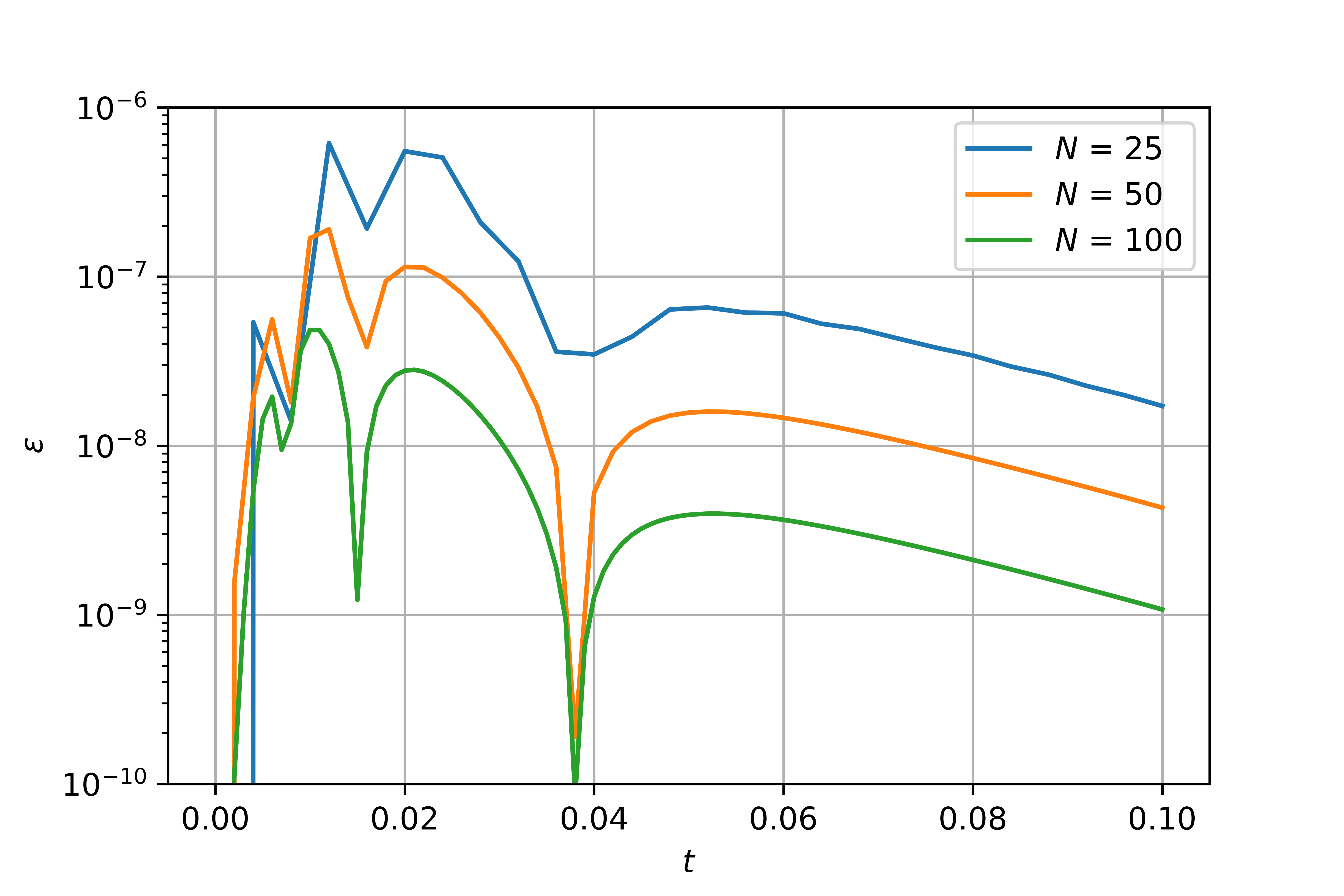} 
\caption{The error of the approximate solution of the problem (\ref{5.1}), (\ref{5.2}) at $\alpha =0.01$ on different time grids.}
\label{f-4}
\end{figure}

The error of the approximate solution of the problem (\ref{5.1}), (\ref{5.2}) was estimated in the grid space norm $L_2(\omega)$ for particular time moments $t = t^n$ as follows:
\[
\varepsilon(t^n) = \Big ( \sum_{x \in \omega} (y^n(x) - u(x,t^n))^2 h \Big )^{1/2} .
\]
The accuracy of the approximate solution using a uniform time grid ($\tau_n = \tau = T/N$) is illustrated by Fig.\ref{f-4}.
We observed the convergence to be approximately second order when a factor of two and four increases the number of steps $N$.

\begin{figure}
\centering
\includegraphics[width=0.75\linewidth]{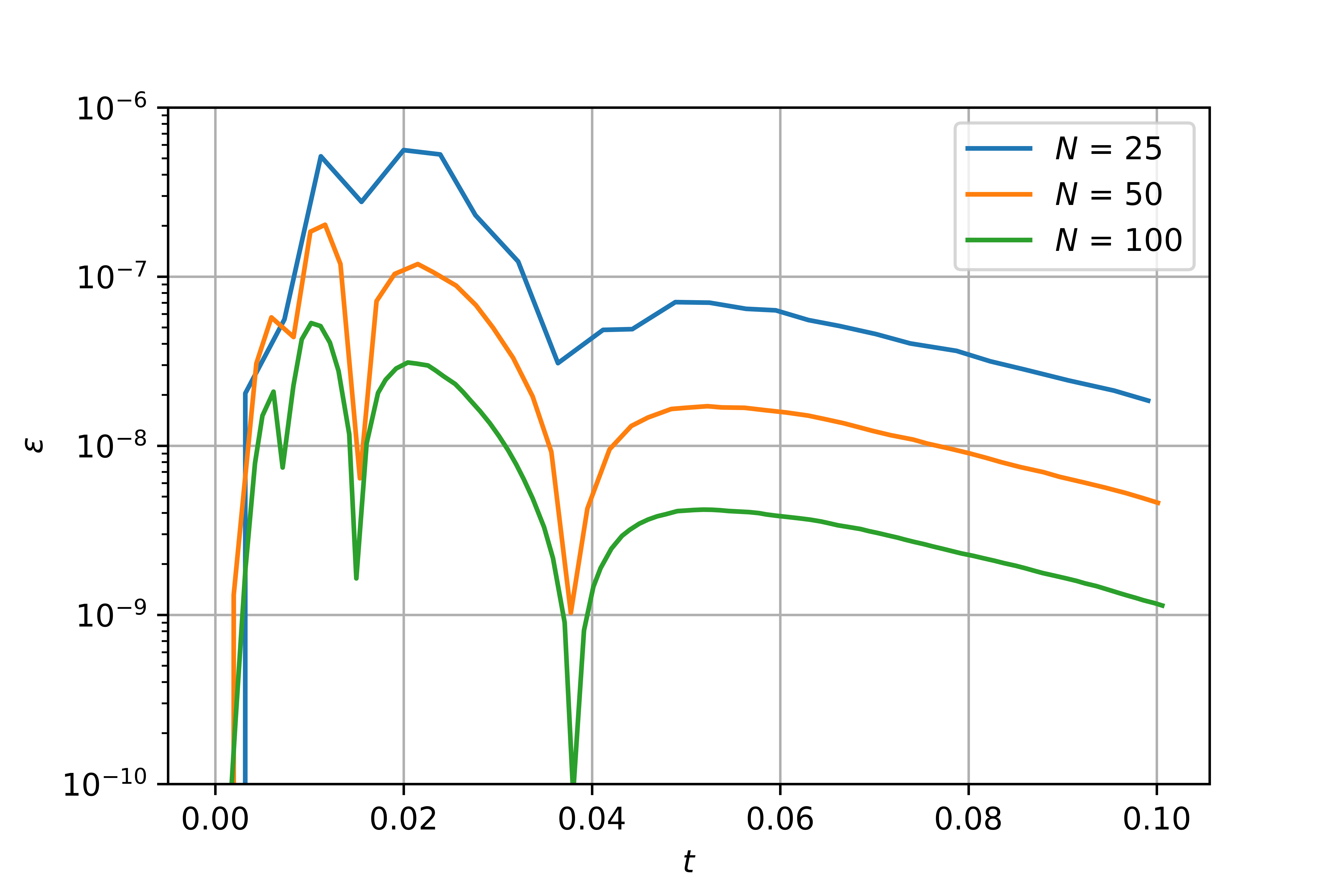} 
\caption{The error of the approximate solution at $\alpha =0.01$ on different time non-uniform grids.}
\label{f-5}
\end{figure}

\begin{figure}
\centering
\includegraphics[width=0.75\linewidth]{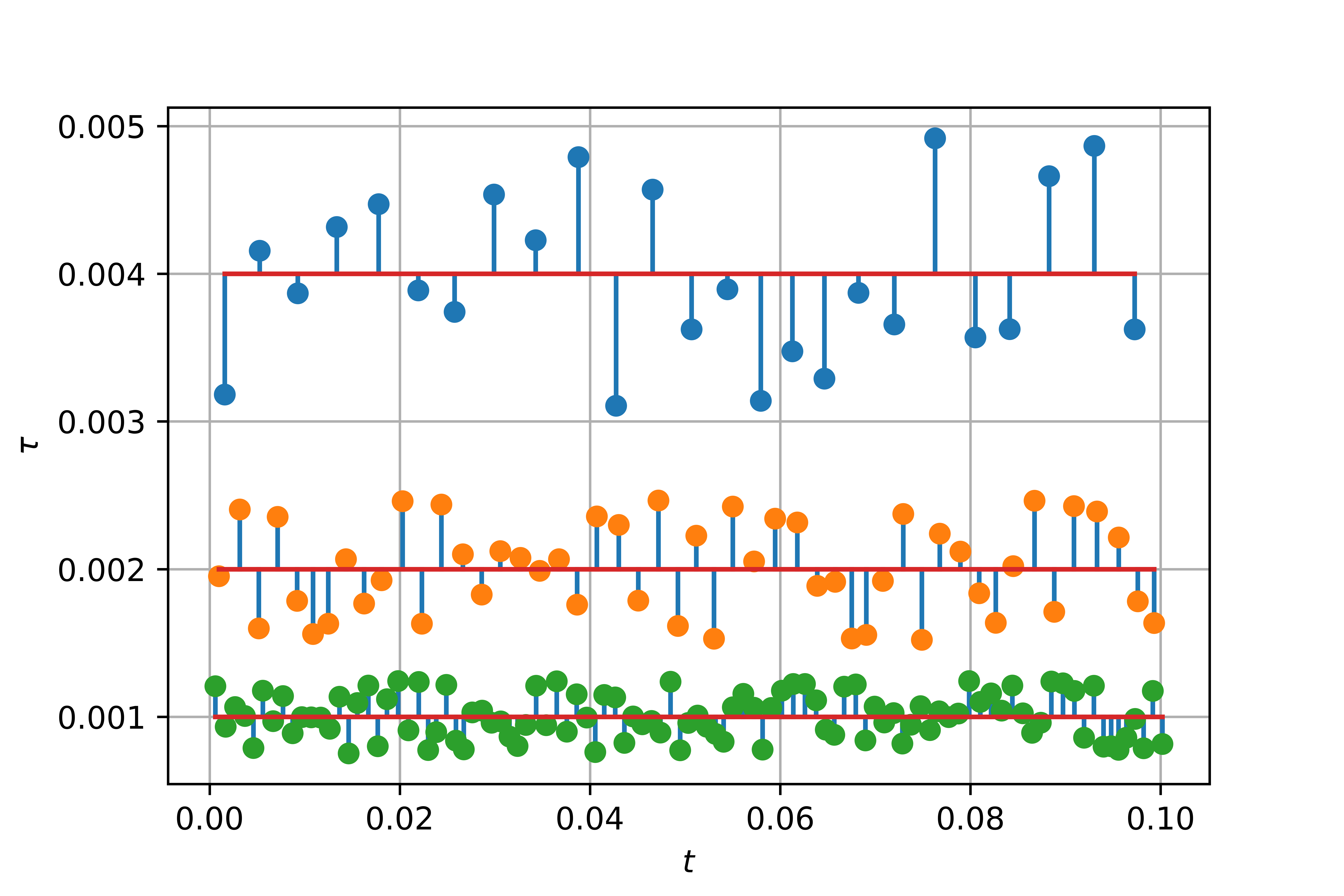} 
\caption{Time grid steps.}
\label{f-6}
\end{figure}

It is convenient to illustrate the capabilities of the constructed schemes on non-uniform time grids with the results obtained by randomly varying the grid step.
The rule calculated the time step
\[
\tau_n = \tau \big (1 + q (\xi - 0.5) \big ),
\quad n = 1,2, \ldots, N ,
\quad \tau = \frac{T}{N} ,
\]
where $\xi$ is a random variable uniformly distributed on the interval $[0,1]$.
In the calculations below, the parameter $q = 0.5$; thus, the grid steps can vary by $5/3$.
The accuracy of the approximate solution using a non-uniform random time grid is shown in Fig.\ref{f-5}.
The non-uniform grids used are shown in Figure \ref{f-6}.
Comparison of Figures \ref{f-4} and \ref{f-5} shows that the use of a substantially non-uniform time grid has little effect on the accuracy of the approximate solution of the considered Cauchy problem for the second-order evolution equation.

\section{Conclusions}\label{sec:6}

\begin{enumerate}
\item
A Cauchy problem is posed for a second-order operator evolution equation in a finite-dimensional Hilbert space with self-adjoint constant operators.
We have well-known results about unconditional stability and convergence of schemes with weights on uniform grids for such problems.
The problem of constructing and investigating three-level operator-difference schemes using non-uniform grids in terms of time is considered.
\item
A second-order evolution equation is written down as a system of two evolution equations with an appropriate choice for new variables.
Two-level schemes with weights on non-uniform time grids are used for approximate solutions to vector problems.
Conditions for unconditional stability of the solution of the Cauchy problem for the system of equations on the initial data and the right-hand side in the corresponding Hilbert spaces have been established.
An equivalent three-level scheme for the original second-order evolution equation has been obtained.
\item
Numerical results of the approximate solution of the Cauchy problem for the evolution equation on non-uniform time grids are presented.
The problem for the bi-parabolic equation, which generalizes classical models of heat conduction, is considered a test one.
The computational data on a grid with randomly changing steps demonstrate the performance of the proposed time approximations.
\end{enumerate}

\section*{Acknowledgements}

This work has been supported by the grants the Russian Science Foundation (RSF 23-41-00037 and RSF 23-71-30013)


\end{document}